\documentclass[11pt]{article} 

\usepackage{pdfsync}
\usepackage[centertags]{amsmath}
\usepackage{amsfonts}
\usepackage{amssymb}
\usepackage{amsthm}
\usepackage{newlfont}
\usepackage[utf8]{inputenc}
\usepackage{mathrsfs}
\usepackage{setspace}
\usepackage{fancyhdr}
\usepackage{epsfig}
\usepackage{graphicx}
\usepackage{accents}

\usepackage{makeidx}
\makeindex

\usepackage[footnotesize,bf]{caption}

\usepackage{verbatim} % for \begin-end{comment} block comments

\DeclareMathSymbol{\widehatsym}{\mathord}{largesymbols}{"62}

\newtheorem{defn}{Definition}[section]

\newtheorem{prop}[defn]{Proposition}

\newtheorem{lem}[defn]{Lemma}

\newcommand{\pf}{\noindent{\bf Proof }\mbox{   }}

\newcommand{\EV}{\mathbb{E}}    % expected value of a random variable
\newcommand{\Z}{\mathbb{Z}}     % integers

\newcommand{\ZK}{\mathbb{Z}^2_K} % 2-dim torus of size K
\newcommand{\bd}{\partial}	    % boundary of a set

\newcommand{\inv}{^{-1}}

\newcommand{\ie}{\emph{i.e.}, }

\newcommand{\st}{such that }

\newcommand{\Dnsc}{D(0,n+s)^c}
\newcommand{\Dnc}{D(0,n)^c}
\newcommand{\Dn}{D(0,n)}

\newcommand{\hp}{\hat{\pi}_K}
\newcommand{\hDnc}{\hp(\Dnc_K)}
\newcommand{\hDnsc}{\hp(\Dnsc_K)}
\newcommand{\hDn}{\hp(\Dn)}

\newcommand{\hdDns}{\hp(\bd D(0,n)_s)}

\newcommand{\hDRc}{\hp(D(0,R)^c_K)}
\newcommand{\hDR}{\hp(D(0,R))}

\newcommand{\hTDnc}{T_{\hDnc}}

\newcommand{\hTdDns}{T_{\hp(\bd D(0,n)_s)}}
\newcommand{\hTDr}{T_{\hp(D(0,r))}}

\newcommand{\hTDRc}{T_{\hp(D(0,R)^c_K)}}

\newcommand{\hGDnc}{\hat{G}_{\hp(\Dnc_K)}}

\newcommand{\hGDn}{\hat{G}_{\hp(D(0,n))}}
\newcommand{\hGDR}{\hat{G}_{\hp(D(0,R))}}

\def\midformat{
\setlength{\itemsep}{0pt} \setlength{\parindent}{0mm}
\setlength{\parskip}{0.12in} \setlength{\textheight}{180mm} % 20110726 to make spaces between paragraphs larger
\setlength{\textwidth}{150mm} \setlength{\evensidemargin}{0in}
\setlength{\oddsidemargin}{0in} \setlength{\topmargin}{0in}
\setlength{\hoffset}{1.0cm} \setlength{\voffset}{0.0cm}
\setlength{\headheight}{15pt} \setlength{\headsep}{.5in}
\setlength{\headwidth}{150mm} } \midformat

\newtheoremstyle{theorem}% name
 {}% Space above
 {}%            Space below
 {\itshape}%    Body font
 {}%            Indent amount (empty = no indent, \parindent = para indent)
 {\ttfamily}%   Thm head font
 {.}%           Punctuation after thm head
 {.5em}%        Space after thm head: " " = normal interword space;
 {}%            Thm head spec (can be left empty, meaning `normal')

\newtheoremstyle{plaintext}% name
 {}% Space above
 {}% Space below
 {\upshape}% Body font
 {}% Indent amount (empty = no indent, \parindent = para indent)
 {\ttfamily}% Thm head font
 {.}% Punctuation after thm head
 {.5em}% Space after thm head: " " = normal interword space;
 {}% Thm head spec (can be left empty, meaning `normal')

\theoremstyle{theorem}

\theoremstyle{plaintext}

\numberwithin{equation}{section}
\begin{document}
\makeatother
%%%%%%%%%%%%%%%%%%%%%%%%%
%   Introduction    %
%%%%%%%%%%%%%%%%%%%%%%%%%
\pagenumbering{roman} \thispagestyle{empty}
\title{Harnack Inequalities of Hitting Distributions \\of Projections of Planar Symmetric Random Walks \\on the Lattice Torus}
\author{Michael Carlisle\footnote{michael.carlisle@baruch.cuny.edu}\\
Baruch College, CUNY}
\maketitle

\begin{abstract}
We give Harnack inequalities for the hitting distributions of a large family of symmetric random walks on $\Z^2$, and their  projections onto the lattice torus $\Z^2_K$. This extends a framework for the simple random walk in \cite{DPRZ2006}, and generalizes the results in \cite{BRFreq} to the toral projection.
\end{abstract}

%-------------

\oddsidemargin 0.0in \textwidth 6.0in \textheight 8.5in

\singlespacing
%\setcounter{tocdepth}{0} % only print ch in toc, no sec or subsec #ing 
%\setcounter{tocdepth}{2} % print ch.sec.subsec in toc
%\tableofcontents 

%\setcounter{tocdepth}{2} % only print ch.sec in toc, no subsec #ing 

%  \cleardoublepage 
 % \addcontentsline{toc}{chapter}{\listfigurename} 
  %\listoffigures 
%\newpage

%  \cleardoublepage 
%\addcontentsline{toc}{chapter}{\indexname}
%\printindex   % why is this breaking ALL formatting??!!

%\doublespacing  %ZZZ UNCOMMENT THIS FOR FINAL VERSION
\pagenumbering{arabic}
\pagestyle{fancy}
\renewcommand{\headrulewidth}{0.0pt}
%\headwidth
\lhead{} \chead{} \rhead{\thepage}
\fancyhead[RO]{\thepage} \fancyfoot{}

\section{Introduction} \label{ch:Intro}

Consider a random walk $S_t = S_0 + \sum_{j=1}^t X_j$, \index{s@$S_t$} on $\Z^2$, with $X = \{X_j\}_{j \in \mathbb{N}}$ \index{x@$X_j$} having the following properties: $X_1$ is symmetric, %recurrent, 
has finite covariance matrix equal to a scalar times the identity, \ie $\Gamma := cov(X_1) = cI$, $c>0$, and $X$ is strongly aperiodic.

As usual in the literature, let 
\begin{align*}
p_1(x,y) = p_1(y-x) = P^x(X_1 = y)
\end{align*}
be the one-step transition probability of $X$. We say $X$ satisfies {\bf Condition A} \index{Condition A} if either $p_1$ has bounded support, or, from any point ``just outside'' a disc, we will enter the disc with positive probability; \ie for any $s \leq n$, for large enough $n$, 
\begin{equation} \label{eqn:ConditionA}
\inf_{y: n \leq |y| < n+s} \sum_{z \in D(x,n)} p_1(y,z)
 = \inf_{y \in \bd D(x,n)_s} P^y(X_1 \in D(x,n)) \geq c e^{-\beta s^{1/4}},
\end{equation}
where the (Euclidean) $s$-annulus around the disc $D(x,n)$ is defined as \index{$\bd D(x,n)_s$}
\begin{equation*} %\label{eqn:s-band}
\bd D(x,n)_s := D(x,n+s) \setminus D(x,n).
\end{equation*}
In particular, if $X_1$ has infinite range, then for any $y \in \bd D(0,n)_s$, there exists $x \in D(0,n)$ such that $p_1(y,x) > 0$.

Starting at a point $x \in A^c$, we define the \emph{hitting distribution} \index{hitting@$H_A$} of $A$ to be 
\[H_A(x,y) := P^x(S_{T_A} = y), \]
where $T_A$ is the first hitting time of the set $A$ by the random walk: 
\[ T_A = \inf\{t \geq 0: S_t \in A \}. \]
The \emph{last exit decomposition} \index{last exit decomposition} of a hitting distribution is based on the Green's function: for $A$ a proper subset of $\Z^2$, $x \in A^c$, and $y \in A$, 
\begin{equation} \label{eq:LastExit}
H_{A}(x,y) = \sum_{z \in A^c} G_{A^c}(x,z) p_1(z,y), 
\end{equation}
where the (truncated) Green's function, up to escaping a set $B$, is defined, for $x, y \in B$, as the total expected number of visits to $y$, starting from $x$, before escaping $B$: 
\begin{equation} \label{eq:GreenDefn}
G_B(x,y) := \EV^x\bigg[\sum_{j=0}^{\infty} 1_{\{S_j=y; j<T_{B^c}\}}\bigg] = \sum_{j=0}^{\infty} P^x(S_j = y; j < T_{B^c})
\end{equation}
and 0 if $x$ or $y \not \in B$. A useful identity relates the Green's function of a set to its expected escape time: if $x \in B$, 
\begin{equation} \label{eq:GreenHit}
\EV^x(T_{B^c}) = \sum_{z \in B} G_B(x,z).
\end{equation}

The goal of this paper is to establish Harnack inequalities of hitting distributions of discs and disc complements on the planar and toral lattices.
%, effectively stating that, up to a small error, the probabilities of having a particular endpoint of a path when starting far inside a disc (if tracking disc escape) or far outside a disc (if tracking disc entry) are nearly the same, making the hitting distribution nearly uniform for large radii. 
Condition A is sufficient to allow the error terms in our Harnack inequalities induced by the toral projection to drop out of sight, as long as there is a moment condition in place on the size of our jumps: we assume 
\begin{equation} \label{eqn:moments}
\EV|X_1|^{M} = \sum_{x \in \Z^2} |x|^{M} p_1(x) < \infty
\end{equation}
for some $M > 4$ (and write $M = 4 + 2 \beta$ for some $\beta > 0$). While $M > 2$ suffices for interior Harnack inequality results, $M = 3 + 2\beta$ is needed for results on the plane, and one more moment is used in our arguments in the transition to the torus. 

%REDUNDANT
%The arguments presented here display a heavy interplay between the notions of visiting a point on the planar lattice $\Z^2$ and examining its projection onto the lattice torus $\Z^2_K$; the Markov-inequality-based error involved in so-called ``targeted jumps'' which, for example, in the projection from $\Z^2 \to \Z^2_K$, land a planar-escape jump back into the toral disc, are dampened sufficiently by this choice of finite moments.

We will switch between the planar and toral representations of the random walk and corresponding stopping times, hitting distributions, etc. %Hats will denote the versions of variables on $\Z_K^2$; no hat means the object is planar. 
Define the projections, \index{pi@$\pi_K$} \index{pihat@$\hp$} for $x = (x_1, x_2) \in \Z^2$, by 
\[\begin{array}{ll}
\pi_K: \Z^2 \to [-K/2,K/2)^2 \cap \Z^2, \\
\pi_K(x) = \left((x_1 + \lfloor \frac{K}{2} \rfloor) (\text{mod } K) - \lfloor \frac{K}{2} \rfloor, (x_2 + \lfloor \frac{K}{2} \rfloor) (\text{mod } K) - \lfloor \frac{K}{2} \rfloor\right); \\
 \hp: \Z^2 \to \Z^2_K, 
\,\,\,\, \hp(x) = (\pi_K x) + (K \Z)^2.
\end{array}\]
(For example, if $x = (-12,6)$ and $K = 11$, then $\pi_{11}(\Z^2) = \{-5,\ldots,5\}^2$, $\pi_{11}(x) = (-1,-5)$, and $\hat \pi_{11}(x) = (-1,-5) + (11\Z)^2$.)

We call the set of lattice points $\pi_K(\Z^2) = [-K/2,K/2)^2 \cap \Z^2$ the {\bf primary copy} \index{primary copy} in $\Z^2$,
% any $x$ in the primary copy  the {\bf primary} $x \in \Z^2$,
 and for $x \in \pi_K(\Z^2)$, $\hat{x} := \hp x$ is its corresponding element in $\Z^2_K$. Any $z \in \pi_K \inv x$, $z \neq \pi_K x$, is called a {\bf copy} \index{copy} of $x$. Likewise, for a set $A \subset \Z^2$, $\hat{A} := \hp A$ is the toral projection of $A$, and the set of all copies of $A$ is \index{piinv@$\hp\inv \hat{A}$}
\[\pi_K \inv\pi_K A = \hp\inv \hat{A} := \{z \in \Z^2: z = x + (iK, jK), \,\,  i, j \in \Z, x \in A\}.\]
%Figure \ref{fig:plane_torus_pullback} displays the projection of a planar set $A$ onto the torus as $\hat{A}$, and its pullback onto $\pi_K \inv A$. (If $A \subset \pi_K \Z^2$, then of course, $A = \pi_K A$.)

For a given $\hat{x} \in \ZK$, we define $x$ to be the (planar) primary copy of that element; $x := \pi_K \hp \inv \hat x$.

While $X_j$ is the $j$th step of the planar walk and $S_j$ its position at time $j$, we use $\hat{S}_j := \hp S_j$ \index{X@$\hat{S}_j$} to denote the position of the toral walk at time $j$. The distance between two points $x,y \in \Z^2$ will be the Euclidean distance $|x-y|$; on the torus, the distance between two points $\hat{x},\hat{y} \in \Z_K^2$ \index{xhat@$\hat{x}$} will be the minimum Euclidean distance $|\hat{x}-\hat{y}| \leq K\sqrt{2}/2$. To limit the issues regarding this distance, we will restrict any discs on $\Z^2_K$ to have radius $n < K/4$ (sometimes written as a diameter constraint: $2n < K/2$).

To bound our functions, we need a precise notion of bounding distance on the lattice torus $\Z^2_K$. As in \cite{DPRZ2006}, a function $f(x)$ is said to be $O(x)$ if $f(x)/x$ is bounded, uniformly in all implicit geometry-related quantities (such as $K$). That is, $f(x) = O(x)$ if there exists a universal constant $C$ (not depending on $K$) such that $|f(x)| \leq Cx$. Thus $x = O(x)$ but $Kx$ is \emph{not} $O(x)$. A similar convention applies to $o(x)$.

Next, we will define a few terms describing the distance of a random walk step, relative to a reference disc of radius $n$ and an $s$-annulus around the disc.
A {\bf small} jump \index{jump}%!small}
 refers to a step that is short enough to possibly (but not necessarily) stay inside a disc of radius $n$ (\ie $|X_1| < 2n$).
A {\bf baby} jump %\index{jump!baby} 
refers to a small jump that is too short to hop over an $s$-annulus from inside a disc (\ie $|X_1| < s$).
A {\bf medium} jump %\index{jump!medium} 
refers to a step that is sufficiently large to hop out of a disc and past an $s$-annulus, but with magnitude strictly less than $K$, and cannot land near a toral copy of its launching point (\ie $s < |X_1| < K-2n$).
A {\bf large} jump %\index{jump!large} 
is a step which, in the toral setting, would be considered ``wrapping around'' in one step (\ie $|X_1| > K-2n$).
A {\bf targeted} jump %\index{jump!targeted} 
is a large jump which lands directly in a copy of the disc or annulus just launched from (\ie $j(K-2n) \leq |X_1| \leq j(K+2n)/\sqrt{2}$ for some $j$). These terms will aid in dealing with differences between planar and toral hitting and escape times.

\section{Random Walk Preliminaries}

In this section we give, without proof, results from \cite{CarVisits} which are used in our Harnack inequality proofs.

First, for $x, y \in \Z^2$ such that $|x| \ll |y|$, we have, by a Taylor expansion around $y$, 
\begin{align} 
\log|y-x| & = \log|y| + O\left( \frac{|x|}{|y|} \right). \label{eq:logEstimate}
\end{align}
In particular, if $x \in D(0,2r)$ and $y \in D(0,R/2)^c$, with $R = 4mr$, we have 
\begin{align} 
\log|y-x| & = \log|y| + O\left( m^{-1} \right). \label{eq:logEstimate2}
\end{align}
Note that \eqref{eq:logEstimate} and \eqref{eq:logEstimate2} hold in the toral case without adjustment.

\subsection{Expected Hitting Times}

Here we state some results about the expected escape and entry times of discs on the plane and torus. Our first generalizes a common disc escape argument, and improves on \cite[Prop. 6.2.6]{LawSRW}.

\begin{lem} \label{lem:EscapeDisc} \cite[Lemma 2.1]{CarVisits}
Let $S_t = S_0 + \sum_{j=1}^t X_j$ be a random walk in $\Z^2$ with $E|X_1|^2 < \infty$, and covariance matrix $\Gamma$ such that $tr(\Gamma) = \gamma^2 > 0$. Then,
uniformly for $x \in D(0,n)$, and for sufficiently large $n$, \index{expected@$\EV^x(T_{D(0,n)^c})$}
\begin{equation} \label{eqn:EscapeDiscExp}
\frac{n^2 - |x|^2}{\gamma^2} \leq \EV^x(T_{D(0,n)^c}) \leq \frac{n^2 - |x|^2}{\gamma^2} + 2n + 1.
\end{equation}
\end{lem} % also, now see \cite[Prop. 6.2.6, p. 152]{LawSRW}

Computational bounds on the $\EV^{\hat x}(\hTDnc)$, the expected toral disc escape, have a slight torally-induced error term \cite[(2.29)]{CarVisits}: 
\begin{equation} \label{eqn:EscapeDiscExpTorus}
\frac{n^2 - |x|^2}{\gamma^2} \leq \EV^{\hat x}(\hTDnc) \leq \frac{n^2 - |x|^2}{\gamma^2} + 2n + 1 + O(K^{-M} n^4).
\end{equation}

\subsection{Probability Estimates}

By Markov's inequality, large jumps are rare: if $C_M = \EV(|X_1|^M) < \infty$, then since $2n < K/2$, 
\begin{equation} \label{eqn:LargeJumpEstimate}
P(|X_1|>K-2n) \leq \frac{C_M}{(K-2n)^M} < \frac{2^M C_M}{K^M} = O(K^{-M}).
\end{equation}
This leads to the targeted jump probability estimate, \ie the rare chance that a large jump lands back into the toral disc the walk escaped from in the planar setting \cite[(2.24)]{CarVisits}: 
\begin{align}
P^x(\hTDnc>T_{D(0,n)^c}) & = \sum_{z \in \left(\hp\inv\hDn \,\setminus\, D(0,n)\right)} \sum_{y \in D(0,n)} G_{D(0,n)}(x,y) p_1(y,z) \notag\\
 & \leq c K^{-M} \sum_{y \in D(0,n)} G_{D(0,n)}(x,y) = O(K^{-M} n^2). \label{eqn:ProbTorusExit}
\end{align}
\cite[(2.50)]{CarVisits} gives a toral gambler's ruin probability estimate for a radius-ruin, \ie hitting the center of a disc before escaping it, in $\Z^2_K$: 
\begin{align}
P^{\hat x}(T_{\hat{0}} < \hTDnc) & = \frac{\log(n/|\hat{x}|) + O(|\hat{x}|^{-1/4})}{\log(n)}\bigg(1 + O((\log n)^{-1})\bigg)   + O(K^{-M} n^2) \notag\\
 & = \frac{\log(n/|\hat{x}|) + O(|\hat{x}|^{-1/4})}{\log(n)}\bigg(1 + O((\log n)^{-1})\bigg). \label{eq:ProbZeroBeforeDisc}
\end{align}
Arguments similar to the proof of \eqref{eq:ProbZeroBeforeDisc} also result in probabilities of ``near ruin'', \ie entering a smaller disc (rather than hitting the center), and their ``success'' counterparts, in both $\Z^2$ and $\Z^2_K$: from \cite[(2.51)-(2.54)]{CarVisits}, uniformly for $r < |x| < R$,  
\begin{align} 
P^x(T_{D(0,r)} > T_{D(0,R)^c}) & = \frac{\log(|x|/r) + O(r^{-1/4})}{\log(R/r)} \label{eq:FP(2.20)} \\
P^x(T_{D(0,r)} < T_{D(0,R)^c}) & = \frac{\log(R/|x|) + O(r^{-1/4})}{\log(R/r)} \label{eq:FP(2.21)} 
\end{align}
\begin{align} 
P^{\hat x}(\hTDr > \hTDRc) & = \frac{\log(|\hat{x}|/r) + O(r^{-1/4})}{\log(R/r)} + O(K^{-M} R^2) \notag\\
 & = \frac{\log(|\hat{x}|/r) + O(r^{-1/4})}{\log(R/r)} \label{eq:ToralGamblersSuccess}  \\
P^{\hat x}(\hTDr < \hTDRc) & = \frac{\log(R/|\hat{x}|) + O(r^{-1/4})}{\log(R/r)} + O(K^{-M} R^2) \notag\\
 & = \frac{\log(R/|\hat{x}|) + O(r^{-1/4})}{\log(R/r)}. \label{eq:ToralGamblersRuin}
\end{align}

Next, we give bounds from \cite{CarVisits} on the probabilities of exiting a disc when starting far inside, or entering a disc from far outside it, by jumping over an annulus. % ZZZ PREVIOUSLY \eqref{eq:sigmax}. 

\begin{lem} \label{lem:BdEscapeEstimateToral} \cite[Lemmas 4.1, 4.2]{CarVisits}
For $n$ sufficiently large, 
\begin{align} 
%\sup_{x \in D(0,n/2)} & P^x(T_{\bd D(0,n)_s} > T_{D(0,n+s)^c}) & \leq c(s^{-M+2} \lor n^{-M+2}). \label{eq:BdEscapeEstimatePlanar} \\
\sup_{\hat x \in \hp(D(0,n/2))} & P^{\hat x}(\hTdDns > T_{\hDnsc}) \leq c (s^{-M+2} \lor n^{-M+2}) \label{eq:BdEscapeEstimateToral} \\
\sup_{\hat x \in \hDnsc} & P^{\hat x}(T_{\hDn} < T_{\hdDns}) \leq cn^2 \log(n)^2 (s^{-M} + n^{-M}). \label{eq:FP(2.71)}
\end{align}
\end{lem}

\subsection{Green's Functions}

To calculate Green's functions in $\Z^2$, we require the \emph{potential kernel} of $X$: for our walks in $\Z^2$, this function is defined by, for $x \in \Z^2$, 
\begin{align}
a(x) := \lim_{n \to \infty} \sum_{j=0}^n [p_j(0) - p_j(x)]. \label{eq:axyDefn}
\end{align}
For our class of random walks, %in \cite[(4.28)]{LawSRW}
 the potential kernel can be shown to be 
\begin{align}
a(x) = \frac{2}{\pi_{\Gamma}} \log|x| + C(p_1) + o(|x|^{-1}) \label{eq:axyCIgen}
\end{align}
where $C(p_1)$ is a constant depending on $p_1$ but not $x$, and $\pi_{\Gamma} = 2\pi\sqrt{\det \Gamma}$.

\cite[(2.45)-(2.46)]{CarVisits} give computational results for the \emph{internal} Green's function at zero on $\Z^2$ and $\Z^2_K$ before escaping a disc: 
\begin{align}
G_{D(0,n)}(0,0) & = \frac{2}{\pi_{\Gamma}}\log n + C' + O(n^{-1/4}) \label{eq:FP(2.13)} 
\end{align}
\begin{align} 
\hGDn(\hat 0,\hat 0) & = G_{D(0,n)}(0,0)(1 + O(K^{-M} n^2)) \notag\\
 & = \left(\frac{2}{\pi_{\Gamma}}\log n + C' + O(n^{-1/4})\right)(1 + O(K^{-M} n^2)) \notag\\
 & = \frac{2}{\pi_{\Gamma}}\log n + C' + O(n^{-1/4}). \label{eq:GreenZeroTorusVal}
\end{align}

\cite[(2.55)-(2.58)]{CarVisits} give calculations and bounds for $G_{D(0,n)}(x,0)$,  $\hGDn(\hat x,\hat 0)$, $G_{D(0,n)}(x,z)$, and $\hGDn(\hat x,\hat z)$: for $x \in \Dn$ and $\hat x \in \hDn$: for some $C = C(p_1) < \infty$, 
\begin{align} 
G_{D(0,n)}(x,0) & = P^x(T_0 < T_{\Dnc}) \, G_{\Dn}(0,0) \notag\\
 & = \frac{2}{\pi_{\Gamma}}\log \bigg(\frac{n}{|x|}\bigg) + C + O(|x|^{-1/4}), \label{eq:FP(2.18)} \\
\hGDn(\hat x,\hat 0) & = \frac{2}{\pi_{\Gamma}}\log \bigg(\frac{n}{|\hat x|}\bigg) + C + O(|\hat x|^{-1/4}) \label{eq:GreenXZeroTorusCalc} \\
G_{D(0,n)}(x,z) & \leq G_{D(x,2n)}(0,z-x) \leq c \log n. \label{eq:FP(2.19)} \\
\hGDn(\hat x,\hat z) & = G_{D(0,n)}(x,z) + O(K^{-M} n^2 \log n) \leq c \log n.\label{eqn:GreenXZBounds}
\end{align}

\cite[Lemma 2.8]{CarVisits} gives that, for any $0<\delta<\varepsilon<1$, we can find $0<c_{1}<c_{2}<\infty$, such that for  all $\hat x \in \hDn \setminus \hp(D(0,\varepsilon n))$, $\hat y \in \hp(D(0, \delta n))$ and all $n$ sufficiently large such that $2n < K/2$,
\begin{equation} \label{eq:GreenFP2.2G}
c_{ 1}\frac{\rho( \hat x)\lor 1}{n}\leq \hat{G}_{\hDn}(\hat y,\hat x)\leq c_{2}\frac{\rho( \hat x)\lor 1}{n}.
\end{equation}

Finally, for any $x,y \in \pi_K(\Dnc_K)$ \st $|x| \leq |y|$, there are computational bounds on the \emph{external} Green's function before entering a disc: 
\begin{eqnarray} 
G_{D(0,n)^c}(x,y) & \leq c_j \log |x|, \label{eqn:ExtGreenIneq}\\
\hGDnc(\hat x,\hat y) & \leq \hat{c}_j \log |\hat x|, \label{eqn:ExtGreenIneqToral}
\end{eqnarray}
where $c_j, \hat{c}_j$ depend on $j>2$, $c_j \geq \hat{c}_j$, and in the toral case, such that $|\hat x| < (\frac{K}{2})^{1/j}$ (there is no such restriction on the planar case).

\section{Interior Harnack inequalities}

We call our first Harnack inequality ``interior'': the starting points are from the interior of a disc, and we examine the probabilities of escaping from a far larger disc around it. We find the planar version first, then move it to the torus. 

\begin{lem} \label{lem:InteriorHarnack2}
Uniformly for $1 \leq m \ll r$, with $s \ll \frac{r}{4m}$, $x, x' \in D(0,2r)$, $R=4mr$, and $y \in D(0,R)^c$,   
\begin{align} 
H_{D(0,R)^c}(x,y) & = (1 + O(m^{-1})) H_{D(0,R)^c}(x',y) + O(R^{-M} \log R), \label{eq:IntHarn2}
\end{align}
where the error term is completely absorbed, \ie 
\begin{align} 
H_{D(0,R)^c}(x,y) & = (1 + O(m^{-1})) H_{D(0,R)^c}(x',y), \label{eq:IntHarn2NoError}
\end{align}
if $s \leq (\log R)^4$ and $y \in \bd D(0,R)_{s}$. %$y \in \bd D(0,R)_{s} \cup D(0,2R)^c$. %\index{sstar@$s^*$}

Furthermore, if $x \in \bd D(0,r)_{r}$ and $y \in D(0,R)^c$, 
\begin{align} 
 & P^x \left(S_{T_{D(0,R)^c}}=y, \, T_{D(0,R)^c} < T_{D\left(0,\frac{r}{4m}+s\right)}\right) \label{eq:IntHarn2OutFirst} \\
  & \,\,\, = (1 + O(m^{-1})) P^x \left(T_{D(0,R)^c} < T_{D\left(0,\frac{r}{4m}+s\right)}\right) H_{D(0,R)^c}(x,y) + O(R^{-M} \log R), \notag 
\end{align}  %, \, T_{D(0,R)^c} = T_{\bd D(0,R)_{R}} \right) used to be in RHS
with a similar loss of the error term if $y \in \bd D(0,R)_{s}$. 
\end{lem}

\pf (In this proof, we switch freely between $R$ and $4mr$.) First, we decompose $D(0,R)$ and examine $H_{D(0,R)^c}(x,y)$:
\begin{align}
H_{D(0,R)^c}(x,y)  & = \left(\sum_{z \in D(0,2mr)} + \sum_{\stackrel{z \in D(0,3mr)}{\setminus D(0,2mr)}} + \sum_{\stackrel{z \in D(0,4mr)}{\setminus D(0,3mr)}}\right) G_{D(0,R)}(x,z) p_1(z,y). \label{eq:IntHarn2Decomp}
\end{align}

If $X$ is finite range, then for $r$ sufficiently large, the first two sums of \eqref{eq:IntHarn2Decomp} are zero. Otherwise, we bound the Green's function via \eqref{eq:FP(2.18)} and \eqref{eq:FP(2.19)}, and by Markov's inequality, $\sum_{z \in D(0,2mr)} p_1(z,y) \leq c(mr)^{-M} \leq cR^{-M}$. Together, these yield, for some $c < \infty$, 
\begin{align}
 G_{D(0,R)}(x,z) \leq G_{D(0,2R)}(0,z) & \leq c \log R \notag\\
\implies \sum_{z \in D(0,2mr)} G_{D(0,R)}(x,z) p_1(z,y) & \leq c R^{-M} \log R. \label{eq:IntHarn2Sum1Bound}
\end{align}
By \eqref{eq:axyCIgen} and \eqref{eq:logEstimate}, uniformly in $x \in D(0,2r)$ and $y \in D(0,2mr)^c$, 
\begin{align}
a(y-x) & = \frac{2}{\pi_{\Gamma}} \log|y-x| + C' + O(|y-x|^{-1}) \notag\\
 & = \frac{2}{\pi_{\Gamma}} \log|y| + C' + O(m^{-1})  = a(y) + O(m^{-1}). \label{eq:axyHarn2}
\end{align}
For $z \in D(0,4mr) \setminus D(0,2mr)$, by the symmetry of the Green's function, the fact that $H$ is a probability,  \cite[(4.28)]{LawSRW}, and \eqref{eq:axyHarn2}, we have 
% removed \cite[(2.8)]{BRFreq} in favor of \cite[(4.28)]{LawSRW}, \cite[(3.5), and (3.6)]{BRFreq} in favor of \eqref{eq:axyCIgen}
\begin{align}
G_{D(0,R)}(x,z) & = G_{D(0,R)}(z,x) \notag\\
 & = \left( \sum_{w \in D(0,R)^c} H_{D(0,R)^c}(z,w) a(w-x)\right) - a(z-x) \label{eq:GD04mr} \\
 & = \left( \sum_{w \in D(0,R)^c} H_{D(0,R)^c}(z,w) a(w)\right) - a(z) + O(m^{-1})\notag\\
 & = G_{D(0,4mr)}(z,0) + O(m^{-1}). \notag
\end{align} 
By \eqref{eq:FP(2.18)}, $G_{D(0,R)}(z,0) \geq c > 0$ uniformly for $z \in D(0,3mr) \setminus D(0,2mr)$, yielding 
\begin{align}
G_{D(0,R)}(x,z) & = G_{D(0,R)}(0,z) (1 + O(m^{-1})). \label{eq:IntHarn2Sum2Bound}
\end{align} 
For $z \in D(0,4mr) \setminus D(0,3mr)$, by the strong Markov property at $T_{D(0,3mr)}$, 
\begin{align}
G_{D(0,R)}(z,x) & = \EV^z (G_{D(0,R)}(S_{T_{D(0,3mr)}},x) ; T_{D(0,3mr)} < T_{D(0,4mr)^c}) \label{eq:IH2EV3}\\
  & = \EV^z (G_{D(0,R)}(S_{T_{D(0,3mr)}},x) ; T_{D(0,3mr)} < T_{D(0,4mr)^c}; |X_{T_{D(0,3mr)}}| > 2mr) \notag\\
  & \,\, + \EV^z (G_{D(0,R)}(S_{T_{D(0,3mr)}},x) ; T_{D(0,3mr)} < T_{D(0,4mr)^c}; |X_{T_{D(0,3mr)}}| \leq 2mr). \notag
\end{align}
By \eqref{eq:FP(2.19)} and \eqref{eq:FP(2.71)}, the last term here is bounded, for sufficiently large $r$, by %\njc
\begin{align}
c (\log R) & P^z(|X_{T_{D(0,3mr)}}| \leq 2mr) \leq c (\log R) P^z(T_{D(0,2mr)} < T_{\bd D(0,2mr)_{mr}}) \notag\\
 & \leq c (\log R) (2mr)^2 \log(2mr)^2 \big[(mr)^{-M} + (2mr)^{-M}\big] \notag\\
 & \leq c (\log R)^3 R^{-M+2} \leq c R^{-M+2+\beta}. \notag
\end{align}
Applying \eqref{eq:IntHarn2Sum2Bound} to the first term, then switching it back to its original form, yields, for $z \in D(0,4mr) \setminus D(0,3mr)$,
\begin{align}
G_{D(0,R)}(z,x) = (1 + O(m^{-1})) \, G_{D(0,R)}(z,0) + O( R^{-M+2+\beta} ). \label{eq:IntHarn2Sum3BoundA}
\end{align}
The planar version of \eqref{eq:GreenFP2.2G} gives us $G_{D(0,R)}(z,x) \geq \frac{c}{mr}$ for $z \in D(0,4mr) \setminus D(0,3mr)$. This reduces \eqref{eq:IntHarn2Sum3BoundA} to
\begin{align}
G_{D(0,R)}(z,x) = (1 + O(m^{-1})) \, G_{D(0,R)}(z,0). \label{eq:IntHarn2Sum3Bound}
\end{align}
Combining \eqref{eq:IntHarn2Decomp}, \eqref{eq:IntHarn2Sum1Bound}, \eqref{eq:IntHarn2Sum2Bound}, and \eqref{eq:IntHarn2Sum3Bound} yields \eqref{eq:IntHarn2}. 

For \eqref{eq:IntHarn2NoError}, let $y \in \bd D(0,R)_{s}$. The only thing we need to do here is show that the error terms are absorbed, \ie for some $c > 0$, with $M=4+2\beta$, 
\begin{align}
m^{-1} H_{D(0,R)^c}(x,y) & \geq c R^{-M} \log R. \label{eq:PlanarIntHarnErrorAbsorb} %\notag \\
%\implies H_{D(0,R)^c}(x,y) & \geq c n^{-3(M-1-\beta)} e^{-n(M-\beta)} = c n^{-9} e^{-4n}. 
\end{align}
Wlog, we can show this for $x=0$. First note that, for $|z| \leq \frac{R}{100}$, by \eqref{eq:FP(2.18)} we have 
\[ G_{D(0,R)}(z,0) \geq c \log \frac{R}{(R/100)} = c \log 100 \geq c \geq \frac{c}{R} \]
for some $c > 0$, and for $z \in D(0,R) \setminus D(0,R/100)$, by the planar version of \eqref{eq:GreenFP2.2G}, $G_{D(0,R)}(z,0) \geq \frac{c}{R}$ as well. Hence, by this, a last exit decomposition, and \eqref{eqn:ConditionA}, 
\begin{align}
m^{-1} & H_{D(0,R)^c}(0,y) = m^{-1} \sum_{z \in D(0,R)} G_{D(0,R)}(0,z) p_1(z,y) \geq \frac{c}{mR} \sum_{z \in D(0,R)} p_1(z,y) \label{eq:PlanarIntHarnErrorAbsorb}\\
 & \geq \frac{c}{mR} e^{-\beta s^{1/4}} = c (mR)^{-1} e^{-\beta \log R} = c m^{-1} R^{-1-\beta} > c R^{-M} \log R. \notag
\end{align}

%For \eqref{eq:InteriorHarnackZoomIn}, 
To show \eqref{eq:IntHarn2OutFirst}, we start with the decomposition 
\begin{align} 
 & P^x \left(S_{T_{D(0,R)^c}}=y, \, T_{D(0,R)^c} < T_{D\left(0,\frac{r}{4m}+s\right)}\right) \label{eq:PlanarEscapeDistDecomp} \\ % was \label{eq:TorusEscapeDistDecomp} \\
  & \,\,\, = H_{D(0,R)^c}(x,y) - P^x \left(S_{T_{D(0,R)^c}}=y, \, T_{D(0,R)^c} > T_{D\left(0,\frac{r}{4m}+s\right)}\right). \notag
\end{align}
%Note that $S_{T_{D(0,R)^c}}=y \in \bd D(0,R)_{R}$ implies $T_{D\left(0,R\right)^c} = T_{\bd D\left(0,R\right)_{R}}$. 
By the strong Markov property at $T_{D\left(0,\frac{r}{4m}+s\right)}$, %and the fact that hitting $D\left(0,4mr\right)^c$ at $y \in \bd D\left(0,4mr\right)_{4mr}$ implies $T_{D\left(0,4mr\right)^c} = T_{\bd D\left(0,4mr\right)_{4mr}}$, 
\begin{align} 
P^x & \left(S_{T_{D(0,R)^c}}=y, \, T_{D(0,R)^c} > T_{D\left(0,\frac{r}{4m}+s\right)}\right) \notag \\
%P^x \left(S_{T_{D(0,R)^c}}=y, \, T_{D(0,R)^c} > T_{D\left(0,\frac{2r}{4m}\right)}\right) \notag\\ %, \, T_{D\left(0,R\right)^c} = T_{\bd D\left(0,R\right)_{R}}\right)
 & = \EV^{x}\bigg[H_{D(0,R)^c}\left(S_{T_{D\left(0,\frac{r}{4m}+s\right)}}, y\right); \, T_{D(0,R)^c} > T_{D\left(0,\frac{r}{4m}+s\right)}\bigg]. \label{eq:GoInFirstStrongMarkov} %, \, T_{D\left(0,R\right)^c} = T_{\bd D\left(0,R\right)_{R}}
\end{align}
By \eqref{eq:IntHarn2}, uniformly in $w \in D\left(0,2r\right)$, %$w \in D\left(0,\frac{2r}{4m}\right)$,
% {\bf JR: but (\ref{eq:InteriorHarnack}) requires $\hat x, \hat x' \in \hDr$, and here $\hat x \in \hp(\bd D(0,r)_{\sqrt r})$. What you want to say is by (\ref{eq:InteriorHarnack}) with $r$ replaced by $r+\sqrt r$.}
\[H_{D(0,R)^c} \left(S_{T_{D\left(0,\frac{r}{4m}+s\right)}}, y\right) = \left(1 + O\left(m^{-1}\right)\right) H_{D(0,R)^c}(w, y) + O(R^{-M} \log R).\]
By \eqref{eq:FP(2.20)} and \eqref{eq:FP(2.21)}, with $m \gg 1$, uniformly for $x \in \bd D(0,r)_{r}$ (say $|x| = cr$, $1 < c < 2$), $\exists c', c'' > 0$ such that 
\begin{align}
P^x \left(T_{D(0,R)^c} < T_{D\left(0,\frac{r}{4m}+s\right)}\right) & = \frac{\log(c'm) + O\left( \left(\frac{r}{m}\right)^{-1/4} \right)}{\log(c'' m^2)} = \frac{1}{2} + o\left( \left(\frac{r}{m}\right)^{-1/4} \right), \notag\\  %O\left( \left(\frac{r}{m}\right)^{-1/4} (\log m)^{-1} \right), 
P^x \left(T_{D(0,R)^c} > T_{D\left(0,\frac{r}{4m}+s\right)}\right) & = \frac{1}{2} + o\left( \left(\frac{r}{m}\right)^{-1/4} \right),  \label{eq:HalfwayFractalGambler} %O\left( \left(\frac{r}{m}\right)^{-1/4} (\log m)^{-1} \right), 
\end{align}
%so up to a small error, these two can be swapped out for each other. 
so the probabilities are both bounded below by a constant. (The small $m$ case operates similarly, but due to the small constants involved, the lower bound must be reduced; $\frac{1}{4}$ for one of them suffices.) 
Combining these and \eqref{eq:GoInFirstStrongMarkov} into \eqref{eq:PlanarEscapeDistDecomp} yields 
\begin{align*} 
P^x & \left(S_{T_{D(0,R)^c}} = y, \, T_{D(0,R)^c} < T_{D\left(0,\frac{r}{4m}+s\right)}\right) \\
 & = H_{D(0,R)^c}(x, y) - \EV^{x}\bigg[H_{D(0,R)^c}\left(S_{T_{D\left(0,\frac{r}{4m}+s\right)}}, y\right); \, T_{D(0,R)^c} > T_{D\left(0,\frac{r}{4m}+s\right)}\bigg]\\
 & = H_{D(0,R)^c}(x, y) \bigg[ P^x \left(T_{D(0,R)^c} < T_{D\left(0,\frac{r}{4m}+s\right)}\right) + P^x \left(T_{D(0,R)^c} > T_{D\left(0,\frac{r}{4m}+s\right)}\right)\\
 & \,\, - \left(1 + O\left(m^{-1}\right)\right) P^x \left(T_{D(0,R)^c} > T_{D\left(0,\frac{r}{4m}+s\right)}\right) + O(R^{-M} \log R) \bigg] \\
& = H_{D(0,R)^c}(x, y)\left(1 + O\left(m^{-1}\right)\right) P^x \left(T_{D(0,R)^c} < T_{D\left(0,\frac{r}{4m}+s\right)}\right) + O(R^{-M} \log R). \, \qed
\end{align*}
%\njc {\bf  Where in the proof are you using the assumption that $y \in \bd D(0,R)_{R}$?} -- MC response: nowhere. I guess this is a vestige of an earlier version. Generalizing \eqref{eq:IntHarn2OutFirst}...

We now move these results to the torus. 

\begin{prop} \label{lem:InteriorHarnack2Toral}
For large $r$ and $1 \leq m \ll r$ such that $R = 4mr < K/6$ and $s \leq (\log R)^4$, uniformly for $\hat x, \hat x' \in \hp(D(0,2r))$ and $\hat y \in \hDRc$, %{\bf  I am a little worried about using \eqref{eq:TorusPlanarIntHarnResult} here, since \eqref{eq:TorusPlanarIntHarnResult} was proven in a situation where $R$ was much bigger that $r$. You should explain why this doesn't really matter}
\begin{align} 
\hat{H}_{\hDRc}(\hat x,\hat y)
 = & \left(1 + O\left(m^{-1}\right)\right) \hat{H}_{\hDRc}(\hat x',\hat y) \notag\\
 & + O(R^{-M} \log R \lor K^{-M} R^2). \label{eq:IntHarn2Toral}
\end{align}
Furthermore, uniformly in $\hat x \in \hp(\bd D(0,r)_{r})$ and $\hat y \in \hDRc$,%\hp(\bd D(0,R)_{R})$, 
\begin{align}
 P^{\hat x}&(\hat{S}_{\hTDRc} = \hat y, \hTDRc < T_{\hp(D(0,\frac{r}{4m}+s))}) \notag\\
 & = \left(1 + O\left(m^{-1}\right)\right) P^{\hat x}(\hTDRc < T_{\hp(D(0,\frac{r}{4m}+s))}) \hat{H}_{\hDRc}(\hat x,\hat y) \notag\\
 & \,\, + O(R^{-M} \log R \lor K^{-M} R^2). \label{eq:IntHarn2ToralZoomIn}
\end{align}
%{\bf  Note that here you already have $\frac{r}{4m}$!}

If $\hat y \in \hp(\bd D(0,R)_{s})$, %\njc {\bf  say what $s$ is  }
the error term is absorbed in both of these statements. %Also, if $K= O(R)$, then the error term is $O(K^{-M+(2 \lor \beta)})$. 
\end{prop}

\pf As before, wlog, we can take $\hat{x}' = \hat{0}$. Let $s$ be the size of the annulus for $\hat y$. For brevity, set 
\[\begin{array}{ccccccc}
A_p & := & \{S_{T_{D(0,R)^c}} = y\}, & \,\,\, & d_p & := & |S_{T_{D(0,R)^c}}-S_{T_{D(0,R)^c}-1}|,\\
A_t & := & \{\hat{S}_{\hTDRc} = \hat y\}, & \,\,\, & d_t & := & |S_{\hTDRc}-S_{\hTDRc-1}|.\\
\end{array}\]
Note that $x$ and $y$ are the primary copies of $\hat x$ and $\hat y$, and so $|x - y| \leq \frac{K}{\sqrt{2}}$, but that $d_p$ is a planar distance using a planar escape time and $d_t$ is a planar distance using a toral escape time; hence, both can exceed $\frac{K}{\sqrt{2}}$, the maximum distance between two points in $\Z^2_K$.

To prove \eqref{eq:IntHarn2Toral}, first we re-label \eqref{eq:IntHarn2} as 
\begin{equation} \label{eq:PlanarIntHarn}
H_{D(0,R)^c}(x,y) = P^x(A_p) = \left(1 + O\left(\frac{r}{R}\right)\right) P^{x'}(A_p) + O\left(R^{-M} \log R\right).
\end{equation}
We have the decomposition
\begin{align}
P^x(A_p) & = P^x(A_p; d_p < K-2R) + P^x(A_p; d_p \geq K-2R). \label{eq:PlanarIntHarnDecomp}
\end{align}
On the plane, the second term of \eqref{eq:PlanarIntHarnDecomp} is zero for all but the furthest-away $y$ in the primary copy (\ie $K - 2R \leq |x - y| \leq \frac{K}{\sqrt{2}}$); for those $y$, we have, by \eqref{eqn:LargeJumpEstimate}, \eqref{eq:GreenHit}, and \eqref{eqn:EscapeDiscExp}, 
\begin{align}
P^{x}(A_p;\, d_p \geq K-2R) & = \sum_{z \in D(0,R)} G_{D(0,R)}(x, z) P^z(|X_1-z|>K-2R) \notag\\
 & \leq  cK^{-M} R^2. \label{eq:PlanarFarOff}
\end{align}
The toral version can be written using planar distances as a decomposition, but using the toral disc escape time means a further decomposition comparing the planar and toral escape times a la (\ref{eqn:ProbTorusExit}). We decompose $\hat{H}_{\hDRc}(\hat x,\hat y) = P^{\hat x}(A_t)$ as 
\begin{align} 
% & = & & P^{\hat x}(A_t;\, d_t < K-2R) + P^{\hat x}(A_t;\, d_t \geq K-2R)\\
P^{\hat x}(A_t)  & = P^{\hat x}(A_t;\, d_t < K-2R;\, T_{D(0,R)^c} = \hTDRc) \label{eq:TorusIntHarnDecomp}\\
 & \,\, + P^{\hat x}(A_t;\, d_t \geq K-2R;\, T_{D(0,R)^c} = \hTDRc) \notag\\
% & \,\, + P^{\hat x}(A_t;\, d_t < K-2R;\, T_{D(0,R)^c} < \hTDRc) \notag\\
% & \,\, + P^{\hat x}(A_t;\, d_t \geq K-2R;\, T_{D(0,R)^c} < \hTDRc). \notag
 & \,\, + P^{\hat x}(A_t;\, T_{D(0,R)^c} < \hTDRc). \notag
\end{align}
In the torus, the first term of \eqref{eq:TorusIntHarnDecomp} equals the first term of \eqref{eq:PlanarIntHarnDecomp}, plus a large jump error which contains some paths from the second term of \eqref{eq:PlanarIntHarnDecomp} (if $y$ is far): by \eqref{eq:PlanarFarOff},
% --- the walk never leaves the original disc until toral escape via a small jump: 
\begin{align*}
P^{\hat x}(A_t;\,& d_t < K-2R;\, T_{D(0,R)^c} = \hTDRc) \\
% =\mbox{{\bf  Where is the second term of (\ref{eq:PlanarIntHarnDecomp})?}} & H_{D(0,R)^c}(x,y).
 & = P^x(A_p; d_p < K-2R) + P^x(A_p; d_p \geq K-2R; d_t < K-2R) \\
 & = P^x(A_p) + O(K^{-M} R^2).
\end{align*}
The second term of \eqref{eq:TorusIntHarnDecomp} only occurs if the final, escaping jump is large: by %\njc 
\eqref{eqn:LargeJumpEstimate}, \eqref{eq:GreenHit}, and \eqref{eqn:EscapeDiscExpTorus}, just as in \eqref{eq:PlanarFarOff}, 
\begin{align*}
P^{\hat x}(A_t;\, d_t & \geq K-2R;\, T_{D(0,R)^c} = \hTDRc) \leq P^{\hat x}(A_t;\, d_t \geq K-2R) \\
 & = \sum_{\hat{z} \in \hDR} \hGDR(\hat x, \hat z) P^z(|X_1 - z|>K-2R) \leq  cK^{-M} R^2.
\end{align*}
The last term of \eqref{eq:TorusIntHarnDecomp} requires a large jump to have occured. Hence, by \eqref{eqn:ProbTorusExit}, 
\begin{align*}
% & P^{\hat x}(A_t;\, d_t < K-2R;\, T_{D(0,R)^c} < \hTDRc) \\
% & + P^{\hat x}(A_t;\, d_t \geq K-2R;\, T_{D(0,R)^c} < \hTDRc) \\
P^{\hat x}(A_t;\, T_{D(0,R)^c} < \hTDRc) \leq cK^{-M} R^2.
\end{align*}
Therefore, \eqref{eq:TorusIntHarnDecomp} reduces to
\begin{align} 
%\hat{H}_{\hDRc}(\hat x,\hat y) & = H_{D(0,R)^c}(x,y) + O(K^{-M} R^2), \text{ or } \\
P^{\hat x}(A_t) & = P^x(A_p) + O(K^{-M} R^2). \label{eq:TorusPlanarIntHarnResult}
\end{align}
(Due to targeting, this is generalizable to any planar set $B \subset D(0,R)$ for $R < K/4$.)
%This gives one of the bounds in (\ref{eq:HittingProbComp}). 
Combining this with \eqref{eq:PlanarIntHarn} gives us \eqref{eq:IntHarn2Toral}: 
\[\begin{array}{lll}
P^{\hat x}(A_t) & = & P^x(A_p) + O(K^{-M} R^2)\\
 & = & \left(1 + O\left(\frac{r}{R}\right)\right) P^{x'}(A_p) + O(K^{-M} R^2) + O\left(R^{-M} \log R\right)\\
 & = & \left(1 + O\left(\frac{r}{R}\right)\right) P^{\hat x'}(A_t) + O(K^{-M} R^2) + O\left(R^{-M} \log R\right). 
\end{array}\]
The proof of \eqref{eq:IntHarn2ToralZoomIn} follows from the Markov property argument for \eqref{eq:IntHarn2OutFirst}, using the appropriate toral identities: \eqref{eq:IntHarn2Toral} for \eqref{eq:IntHarn2}, and  \eqref{eq:ToralGamblersSuccess}-\eqref{eq:ToralGamblersRuin} for \eqref{eq:FP(2.20)}-\eqref{eq:FP(2.21)}. $\qed$

\section{Exterior Harnack inequality}

%{\bf \njcc throughout - generalizing} 
To aid our construction of an exterior Harnack inequality (moving from outside a disc to a much smaller disc far inside it), we first establish uniform bounds on external Green's functions and probabilities in the torus and plane. Fix $\delta < 1$ and use $r \geq e^n$ for some $n>13$, $R = 4mr$ for some $1 \leq m \ll r$, and $s \leq (\log R)^4$. First, for %$\hat x \in \hp(\bd D(0,R)_{R^{\delta}})$
 $\hat x \in \hp(\bd D(0,R)_{R/100})$ and $\hat y \in \hDRc$, we show that 
\begin{align}
%\hat{G}_{\hp(D(0,r+n^4)^c)}(\hat x,\hat y) \geq c > 0. \label{eq:FP(3.19)}
\hat{G}_{\hp(D(0,r+s)_{K}^c)}(\hat x,\hat y) \geq c > 0. \label{eq:FP(3.19)}
\end{align}
%{\bf   don't you mean $D(0,r+n^4)_{K}^c$? If so, make sure to change throughout.}
Pick some $\hat x_1 \in \hp(\bd D(0,R))$, and, proceeding clockwise, choose points $\hat x_2, ..., \hat x_{36} \in \hp(\bd D(0,R))$ whose rays beginning at $\hat 0$ divide $\hp(\bd D(0,R))$ into 36 approximately equal arcs. The distance between any two adjacent such $\hat x_j$ is, for sufficiently large $R$, approximately $2R \sin(\pi/36) \approx 0.174R$. Thus, using discs of radius $R/5$ (so adjacent circles contain their neighbor's centers), and by \eqref{eq:ToralGamblersRuin} we have for any $j=1,..., 36$
\begin{align}
& \inf_{\hat x \in \hp(D(x_j, R/5))} P^{\hat x}(T_{\hp(D(x_{j+1}, R/5))} < T_{\hp(D(0,r+s))}) \label{eq:FP(3.20)}\\    
  \geq & \inf_{\hat x \in \hp(D(x_{j+1}, 2R/5))} P^{\hat x}(T_{\hp(D(x_{j+1}, R/5))} < T_{\hp(D(\hat x_{j+1},R/2)^c_K)}) \geq c_1 > 0, \notag
\end{align}
%{\bf  please copy (3.20) of FP correctly }
for some $c_1$ independent of $n, r, R$, $n$ large, and $\hat x_{37} = \hat x_1$. 
%{\bf  At the very least, the first inequality requires $r+n^4\leq R/2$. But you haven't specified any relation between these parameters.}

Hence, by the strong Markov property, rotating through the arcs, we have 
\begin{align}
\inf_{j,k} \inf_{\hat x \in \hp(D(x_j, R/5))} P^{\hat x}(T_{\hp(D(x_{k}, R/5))} < T_{\hp(D(0,r+s))}) \geq c_2 := c_1^{36}. \label{eq:FP(3.21)}
\end{align}

\begin{figure}[!ht]
  \centering
    \includegraphics[width=3in]{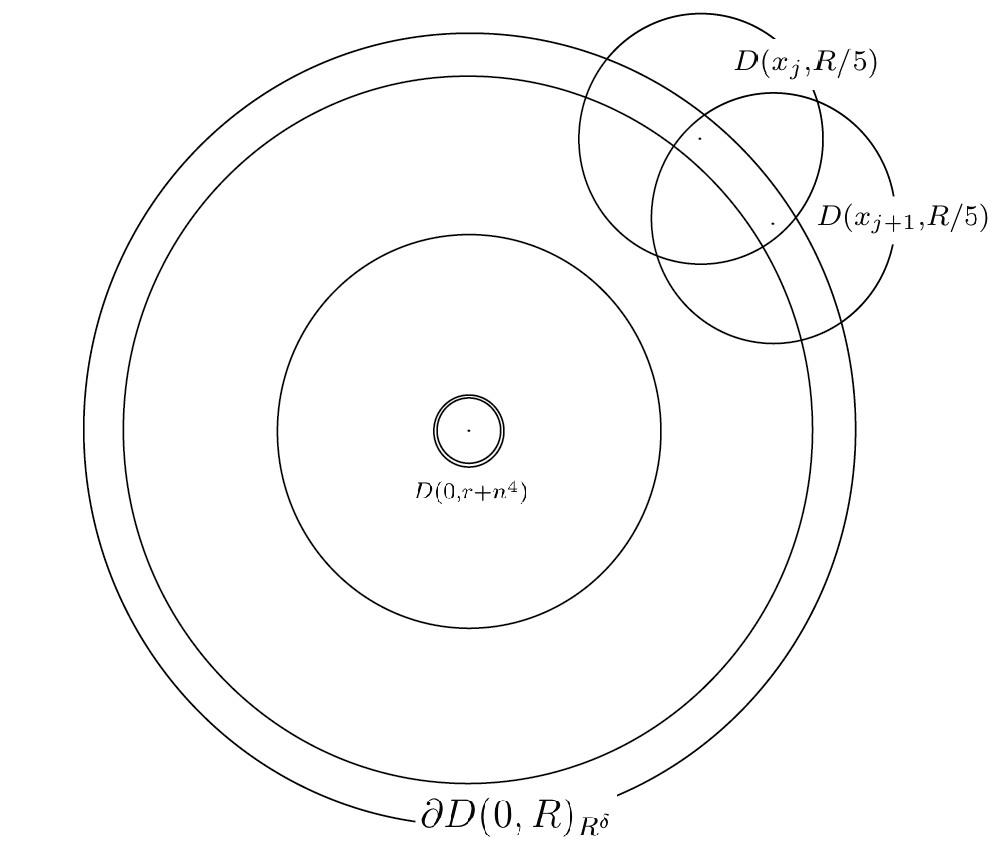}
    \parbox{4in}{
    \caption{Possible selections of $\hat x_j$ and $\hat x_{j+1}$, with associated discs and annuli.}
    \label{fig:Harnack_exterior1}}
\end{figure}

Furthermore, it follows from \eqref{eq:ProbZeroBeforeDisc} that for any $j$, 
\begin{align}
& \inf_{\hat x, \hat x' \in \hp(D(x_j, R/5))} P^{\hat x}(T_{\hat x'} < T_{\hp(D(0,r+s))}) \label{eq:FP(3.22)}\\
 \geq & \inf_{\hat x \in \hp(D(x', 2R/5))} P^{\hat x}(T_{\hat x'} < T_{\hp(D(x',R/2)^c_K)}) \geq \frac{c_3}{\log R} \notag
\end{align}
for some independent $c_3 > 0$. 

Since $\hp(\bd D(0,R)_{R/100}) \subset \cup_{j=1}^{36} \hp(D(x_j, R/5))$, combining \eqref{eq:FP(3.21)} and \eqref{eq:FP(3.22)} we have, for some independent $c_4 > 0$,  
\begin{align}
\inf_{\hat x, \hat x' \in \hp(\bd D(0, R)_{R/100})} P^{\hat x}(T_{\hat x'} < T_{\hp(D(0,r+s))}) \geq \frac{c_4}{\log R}. \label{eq:FP(3.23)}
\end{align}
It then follows from \eqref{eq:GreenZeroTorusVal} that 
\begin{align}
 & \inf_{\hat x, \hat x' \in \hp(\bd D(0, R)_{R/100})} \hat{G}_{\hp(D(0,r+s)^c_K)}(\hat x, \hat x') \notag\\
 = & \inf_{\hat x, \hat x' \in \hp(\bd D(0, R)_{R/100})} P^{\hat x}(T_{\hat x'} < T_{\hp(D(0,r+s))}) \hat{G}_{\hp(D(0,r+s)^c_K)}(\hat x', \hat x') \notag\\
\geq & \frac{c_4}{\log R} \hat{G}_{\hp(D(x',R/2))}(\hat x', \hat x') \geq c_5 > 0 \label{eq:FP(3.24)}
\end{align}
for some independent $c_5 > 0$. Using the strong Markov property, \eqref{eq:FP(3.24)}, and \eqref{eq:FP(2.71)}, we see that 
\begin{align}
& \inf_{\stackrel{\hat z \in \hp(D(0,1.01R)^c_K),}{\hat x \in \hp(\bd D(0,R)_{R/100})}} \hat{G}_{\hp(D(0,r+s)^c_K)}(\hat z,\hat x) \label{eq:FP(3.25)}\\
 & \,\, \geq \EV^{\hat z}\left( \hat{G}_{\hp(D(0,r+s)^c_K)}(\hat{S}_{T_{\hp(D(0,1.01R))}},\hat x); \hat{S}_{T_{\hp(D(0,1.01R))}} \in \hp(\bd D(0,R)_{R/100})\right) \notag\\
 & \,\, \geq c > 0. \notag
\end{align}
This gives \eqref{eq:FP(3.19)}. Applying the same argument once more, 
\begin{align}
& \inf_{\stackrel{\hat z \in \hp(D(0,1.01R)^c_K)}{\hat x \in \hp(D(0,R)^c_K)}} \hat{G}_{\hp(D(0,r+s)^c_K)}(\hat z,\hat x) \label{eq:FP(3.25s)}\\
 & \,\, \geq \EV^{\hat z}\left( \hat{G}_{\hp(D(0,r+s)^c_K)}(\hat{S}_{T_{\hp(D(0,1.01R))}},\hat x); \hat{S}_{T_{\hp(D(0,1.01R))}} \in \hp(\bd D(0,R)_{R/100})\right) \notag\\
 & \,\, \geq c > 0. \notag
\end{align}
Hence, for all $\hat x,\hat y \in \hDRc$, 
\begin{align}
\hat{G}_{\hp(D(0,r+s)_{K}^c)}(\hat x,\hat y) \geq c > 0. \label{eq:FP(3.19s)}
\end{align}

Next, we look at the external Green's function near the $r$-disc: uniformly for $\hat x \in \hp(D(0,R)_{K}^c)$ and $\hat z \in \hp(D(0,2r)) \setminus \hp(D(0,5r/4))$, we have by \eqref{eq:FP(3.19s)} and \eqref{eq:ToralGamblersSuccess}, 
\begin{align}
 & \hat{G}_{\hp(D(0,r+s)^c_K)}(\hat x,\hat z) = \hat{G}_{\hp(D(0,r+s)^c_K)}(\hat z,\hat x) \notag\\
 & \,\, = \EV^{\hat z}\left( \hat{G}_{\hp(D(0,r+s)^c_K)}(\hat{S}_{T_{\hp(D(0,R)^c_K)}},\hat x); T_{\hp(D(0,R)^c_K)} < T_{\hp(D(0,r+s))}\right) \notag\\
 & \,\, \geq c P^{\hat z}(T_{\hp(D(0,R)^c_K)} < T_{\hp(D(0,r+s))}) \geq \frac{c}{\log m}. \label{eq:FP(3.26)}
\end{align}
%{\bf  This uses $\log (R/(r+n^{4}))\leq c \log n$, which is fine but you have to spell out the relation between thes parameters}

Getting closer to the disc, for any $\varepsilon > 0$, uniformly in $\hat x, \hat x' \in \hp(\bd D(0,R)_{R/100})$ and $\hat z \in \hp(D(0,2r)) \setminus \hp(D(0,r + (1 + \varepsilon)s))$, we have by the strong Markov property and \eqref{eq:FP(3.23)}, for any $\hat x'\in \bd D(0,R)_{R/100}$, 
\begin{align}
\hat{G}_{\hp(D(0,r+s)^c_K)}(\hat x,\hat z) & \geq P^{\hat x}(T_{\hat x'} < T_{\hp(D(0,r+s))}) \, \hat{G}_{\hp(D(0,r+s)^c_K)}(\hat x',\hat z) \notag\\
 & \geq \frac{c \hat{G}_{\hp(D(0,r+s)^c_K)}(\hat x',\hat z)}{\log R}. \label{eq:FP(3.28)}
\end{align}
In view of \eqref{eq:GreenFP2.2G}, if $\hat x' \in \hp(\bd D(0,R)_{R/100})$ is chosen as close as possible to the ray from the origin which passes through $\hat z$, we have 
\begin{align}
\hat{G}_{\hp(D(0,r+s)^c_K)}(\hat x',\hat z) & \geq \hat{G}_{\hp(D(\hat x',|\hat x'|-(r+s)))}(\hat x',\hat z) \geq \frac{c}{R}, \label{eq:FP(3.29)}
\end{align}
which, combined with \eqref{eq:FP(3.28)}, gives us 
\begin{align}
 \inf_{\stackrel{\hat x \in \hp(\bd D(0, R)_{R/100})}{\hat z \in \hp(D(0,2r)) \setminus \hp(D(0,r + (1 + \varepsilon)s))}} \hat{G}_{\hp(D(0,r+s)^c_K)}(\hat x,\hat z) & \geq \frac{c}{R \log R}. \label{eq:FP(3.27)}
\end{align}
Using the strong Markov property, \eqref{eq:FP(3.27)}, and \eqref{eq:BdEscapeEstimateToral}, we see that 
\begin{align}
&  \inf_{\stackrel{\hat x \in \hp(D(0,1.01R)^{c}_K))}{\hat z \in \hp(D(0,2r)) \setminus \hp(D(0,r + (1 + \varepsilon)s))}} \hat{G}_{\hp(D(0,r+s)^c_K)}(\hat z,\hat x) \label{eq:FP(3.25t)}\\
 & \,\, \geq \EV^{\hat z}\left( \hat{G}_{\hp(D(0,r+s)^c_K)}(\hat{S}_{T_{\hp(D(0,1.01R)^{c}_K)}},\hat x); \hat{S}_{T_{\hp(D(0,1.01R)^{c}_K)}} \in \hp(\bd D(0,R)_{R/100})\right) \notag\\
 & \,\, \geq \frac{c}{R \log R}. \notag
\end{align}
Hence  
\begin{align}
 \inf_{\stackrel{\hat x \in \hp(D(0,R)^{c}_K))}{\hat z \in \hp(D(0,2r)) \setminus \hp(D(0,r + (1 + \varepsilon)s))}} \hat{G}_{\hp(D(0,r+s)_{K}^c)}(\hat x,\hat y) \geq \frac{c}{R \log R}. \label{eq:FP(3.19t)}
\end{align}

By removing the (hidden) targeted jump error terms, the entire argument in \eqref{eq:FP(3.19)}-\eqref{eq:FP(3.19t)} also applies to the plane. We now find a general planar Harnack inequality for entering a small disc from far outside. 

\begin{prop} \label{lem:ExteriorHarnackPlanarGeneral}
%Let $1/2 \leq \delta < 1$, 
Let $R = 4mr$ with $1 \leq m \ll r$ ($m = o(r^{1/4})$) and large enough $r$, and $s \leq (\log R)^4$. Then, uniformly for $x, x' \in D(0,R)^c$ and $y \in \bd D(0,r)_{s}$, 
\begin{align} 
H_{D(0,r+s)}(x,y) & = \left(1 + O\left(m^{-1} \log m\right)\right) H_{D(0,r+s)}(x',y). \label{eq:ExteriorHarnackPlanarGeneral}
\end{align}
Furthermore, for $x, x' \in \bd D(0,R)_{\sqrt{R}}$, 
\begin{align}
P^x& (S_{T_{D(0,r+s)}} = y; \, T_{D(0,r+s)} < T_{{D(0,4mR)}^c}) \label{eq:ExteriorHarnackPlanarZoomOutGeneral} \\
&  = \left(1 + O\left(m^{-1} \log m\right)\right) H_{D(0,r+s)}(x,y) P^x(T_{D(0,r+s)} < T_{{D(0,4mR)}^c})  \notag\\
&  = \left(1 + O\left(m^{-1} \log m\right)\right) P^{x'}(S_{T_{D(0,r+s)}} = y; T_{D(0,r+s)} < T_{{D(0,4mR)}^c}). \notag
\end{align}
\end{prop}

\pf % We will first prove \eqref{eq:ExteriorHarnackPlanarGeneral} for \njc for $x,x' \in \bd D(0,R)_{R^{\delta}}$
For $x,x' \in D(0,R)^c$ and $y \in \bd D(0,r)_{s}$, we have the last exit decomposition 
\begin{align}
H_{D(0,r+s)}(x,y) & = \left( \sum_{\stackrel{z \in D(0,5r/4)}{\setminus D(0,r+s)}} + \sum_{\stackrel{z \in D(0,2r)}{\setminus D(0,5r/4)}} + \sum_{z \in D(0,2r)^c} \right) G_{D(0,r+s)^c}(x,z) p_1(z,y). \label{eq:FP(3.33)}
\end{align}
Let $x, x' \in \bd D(0,R)_{R}$ and set $N \geq 4mR$. Uniformly for $z \in D(0,2r) \cup D(0,N)^c$, by \eqref{eq:axyCIgen} and \eqref{eq:logEstimate2},  
\begin{align}
 a(x-z) = \frac{2}{\pi_{\Gamma}} \log|x-z| + C' + O(|x-z|^{-1}) = a(x'-z) + O(m^{-1}). \label{eq:axzR}
\end{align}
Using the same approach as in \eqref{eq:GD04mr}, \eqref{eq:axzR} implies that, %\cite[(3.36)]{BRFreq} holds: 
for $A(r+s,N) := D(0,N) \setminus D(0,r+s)$, 
\begin{align*}
G_{A(r+s,N)}(x,z) = G_{A(r+s,N)}(x',z) + O(m^{-1}),  
\end{align*}
which, by letting $N \to \infty$ and applying the dominated convergence theorem, 
\begin{align}
G_{D(0,r+s)^c}(x,z) = G_{D(0,r+s)^c}(x',z) + O(m^{-1}). \label{eq:FP(3.37)}
\end{align}
Applying \eqref{eq:FP(3.26)} to \eqref{eq:FP(3.37)} yields, for $z \in D(0,2r) \setminus D(0,5r/4)$, % ZZZ I just want to note here that, as in the interior Harnack inequality proof in about the same spot, this is big-O here, but it should probably be little-o, since the error is just O(n^{-3}). We only have a LOWER bound of c/\log n \leq G_{D(0,r+n^4)^c}(x',z), so we probably need to shrink G_{D(0,r+n^4)^c}(x',z) FASTER than n^{-3} \log n -- G_{D(0,r+n^4)^c}(x',z) could be as big as its upper bound of c \log r \geq cn, a LOT bigger than c/\log n. 
\begin{align}
G_{D(0,r+s)^c}(x,z) = (1 + O(m^{-1} \log m)) G_{D(0,r+s)^c}(x',z). \label{eq:FP(3.34)}
\end{align}
Next, by the symmetry of the Green's function, the strong Markov property at $T_{D(0,5r/4)^c}$, \eqref{eq:FP(3.34)} for $z \in D(0,5r/4) \setminus D(0,r+s)$, and decomposing, we have 
\begin{align}
 & G_{D(0,r+s)^c}(x,z) = G_{D(0,r+s)^c}(z,x) \label{eq:FP(3.39)}\\
 & = \EV^z( G_{D(0,r+s)^c}(S_{T_{D(0,5r/4)^c}},x) \, ; \, T_{D(0,5r/4)^c} < T_{D(0,r+s)} ), \notag \\
 & = \EV^z( G_{D(0,r+s)^c}(S_{T_{D(0,5r/4)^c}},x) \, ; \, T_{D(0,5r/4)^c} < T_{D(0,r+s)} \, , \, |S_{T_{D(0,5r/4)^c}}| \leq 2r)  \notag \\
 & \,\,\, + \, \EV^z( G_{D(0,r+s)^c}(S_{T_{D(0,5r/4)^c}},x) \, ; \, T_{D(0,5r/4)^c} < T_{D(0,r+s)} \, , \, |S_{T_{D(0,5r/4)^c}}| > 2r).  \notag
\end{align}
By \eqref{eq:FP(3.34)} on the first term and \eqref{eqn:ExtGreenIneq} on the second term, \eqref{eq:FP(3.39)} is bounded above: 
%{\bf  but that requires $r<K^{1/3}$} -- no, this is planar; there are no restrictions on the planar bound (3.9)}, by {\bf  why a bound? Is this an equality?}
\begin{align}
G_{D(0,r+s)^c}(z,x) \leq & (1 + O(m^{-1} \log m)) \label{eq:FP(3.39)-(3.40)} \\
 & \,\, \EV^z( G_{D(0,r+s)^c}(S_{T_{D(0,5r/4)^c}},x') \, ; \, T_{D(0,5r/4)^c} < T_{D(0,r+s)} ) \notag\\
 & \,\,\,\, + c \log(R) P^z( |S_{T_{D(0,5r/4)^c}}| > 2r ). \notag
\end{align}
Applying the first two lines of \eqref{eq:FP(3.39)} again, the first term here is 
\begin{align*}
(1 &+ O(m^{-1} \log m)) \EV^z( G_{D(0,r+s)^c}(S_{T_{D(0,5r/4)^c}},x') \, ; \, T_{D(0,5r/4)^c} < T_{D(0,r+s)} ) \\  
& = (1 + O(m^{-1} \log m))G_{D(0,r+s)^c}(z,x') = (1 + O(m^{-1} \log m))G_{D(0,r+s)^c}(x',z).
\end{align*}
A last exit decomposition of $P^z( |S_{T_{D(0,5r/4)^c}}| > 2r )$, then \eqref{eq:GreenXZeroTorusCalc} and \eqref{eqn:moments} yield 
\begin{align*}
G_{D(0,r+s)^c}(x,z) \leq & (1 + O(m^{-1} \log m)) G_{D(0,r+s)^c}(x',z) \\
 & \,\,\, + c \log(R) \sum_{\stackrel{|y|<5r/4}{2r<|w|}} G_{D(0,5r/4)}(z,y) p_1(y,w) \\
  \leq & (1 + O(m^{-1} \log m)) G_{D(0,r+s)^c}(x',z) + c \log(R) \log(r) r^{-M+2}.
\end{align*}
Since this argument is symmetric in $x$ and $x'$, then we have that for $z \in D(0,5r/4) \setminus D(0,r+s)$, and $c \log R = c(\log 4 + \log m + \log r) = O(\log r)$, 
\begin{align}
G_{D(0,r+s)^c}(x,z) = (1 + O(m^{-1} \log m)) G_{D(0,r+s)^c}(x',z) + O(r^{-M+2} (\log r)^2). \label{eq:FP(3.38)}
\end{align}
Finally, by \eqref{eqn:ExtGreenIneq}, %{\bf  same problem -- again, it's planar, not a problem here} 
for $z \in D(0,2r)^c$, %by \eqref{eq:FP(2.71)} {\bf  where is this used?  }, 
$G_{D(0,r+s)^c}(x,z) = O(\log R)$. Thus, for $y \in \bd D(0,r)_{s}$, since $\sum_{z \in D(0,2r)^c} p_1(z,y) \leq O(r^{-M})$ by symmetry and Markov's inequality, 
\begin{align}
\sum_{z \in D(0,2r)^c} G_{D(0,r+s)^c}(x,z) p_1(z,y) = O(r^{-M} \log r). \label{eq:FP(3.41)}
\end{align}
Combining \eqref{eq:FP(3.34)}, \eqref{eq:FP(3.38)}, and \eqref{eq:FP(3.41)} bounds the sums in \eqref{eq:FP(3.33)} to 
\begin{align}
H_{D(0,r+s)}(x,y) & = (1 + O(m^{-1} \log m)) H_{D(0,r+s)}(x',y) + O(r^{-M+2} (\log r)^2). \label{eq:ExtHarnAlmost}
\end{align}

To complete the proof of \eqref{eq:ExteriorHarnackPlanarGeneral} for $x,x' \in \bd D(0,R)_{R}$, we must show that, uniformly for $x \in \bd D(0,R)_{R}$ and $y \in \bd D(0,r)_{s}$, 
\begin{align}
r^{-M+2} (\log r)^2 \leq c(m^{-1} \log m) H_{D(0,r+s)}(x,y). \label{eq:FP(3.42)}
\end{align}
With $A_r := D(0,2r) \setminus D(0,r+(1+\varepsilon)s)$, using a last exit decomposition 
%, applying the strong Markov property at $T_{D(0,R+R^{\delta})}$, 
and bounding with the planar version of \eqref{eq:FP(3.19t)},
% and then \eqref{eq:FP(2.71)} yields 
\begin{align}
H_{D(0,r+s)}(x,y) & = \sum_{z \in D(0,r+s)^c} G_{D(0,r+s)^c}(x,z) p_1(z,y) \label{eq:FP(3.43)}\\
 & \geq \sum_{z \in A_r} G_{D(0,r+s)^c}(x,z) p_1(z,y) 
% & \geq \sum_{z \in A_r} E^{x}\left[G_{D(0,r+s)^c}\left(S_{T_{D(0,R+R^{\delta})}}, \, z\right)\right] p_1(z,y) \notag \\
% & \geq \sum_{z \in A_r} E^{x}\left[G_{D(0,r+s)^c}\left(S_{T_{\bd D(0,R)_{R^{\delta}}}}, \, z\right); \, T_{D(0,R+R^{\delta})} = T_{\bd D(0,R)_{R^{\delta}}}\right] p_1(z,y) \notag \\
% & \geq \frac{cP^{x}(T_{D(0,R+R^{\delta})} = T_{\bd D(0,R)_{R^{\delta}}})}{R \log R} \sum_{z \in A_r}  p_1(z,y) \notag \\
% & \geq \frac{c(1 - c'R^{-M\delta+2}(\log R)^2)}{R \log R} \sum_{z \in A_r}  p_1(z,y)
 \geq \frac{c''}{R \log R} \sum_{z \in A_r}  p_1(z,y) \notag
\end{align}
for any $\varepsilon > 0$. 
Note that the annulus $A_r$ contains the disc $D(v,2(1 + \varepsilon)s)$, where $v:= (r+3(1+\varepsilon)s)y/|y|$. Thus, $2(1 + \varepsilon)s \leq |y-v| \leq 3(1+\varepsilon)s$, and \eqref{eqn:ConditionA} (where we consider $y \in \bd D(v,2(1 + \varepsilon)s)_{(1+\varepsilon)s}$), and with $s \leq (\log R)^4 \leq c(\log r)^4$, 
\begin{align}
\sum_{z \in A_r} p_1(z,y) & \geq \sum_{z \in D(v,2(1 + \varepsilon)s)} p_1(z,y) \geq ce^{-\beta((1+\varepsilon)s)^{1/4}} \geq cr^{-(1+\varepsilon)^{1/4} \beta}. \label{eq:smallsumForCondA}
\end{align}
Hence, combining \eqref{eq:FP(3.43)} and \eqref{eq:smallsumForCondA}, and since $m \leq \sqrt{r}$, $R < r^2$, some $\varepsilon' > 0$, and $r^{\beta} > (\log r)^3$ for large enough $r$, 
%since $r \geq e^n$ and $r^2 > R = n^3 r$, for large enough $n$, 
\begin{align}
c(m^{-1} \log m) H_{D(0,r+s)}(x,y) & \geq \frac{c(m^{-1} \log m) r^{-(1+\varepsilon)^{1/4} \beta}}{R \log R} \notag\\
 & \geq cr^{-(1+\varepsilon)^{1/4} \beta} (m^{-1} \log m) (mr)^{-1} (2 \log r)^{-1} \notag\\
 & \geq cr^{-1-(1+\varepsilon')\beta} m^{-2} (\log m)(\log r)^{-1} \notag\\
 & \geq cr^{-2-2\beta} (\log m)(\log r)^2 \geq cr^{-M+2}(\log r)^2, \label{eq:ExteriorErrorBound}
\end{align}
which proves \eqref{eq:FP(3.42)}, and hence \eqref{eq:ExteriorHarnackPlanarGeneral}, for $x, x' \in \bd D(0,R)_{R}$.

Next we show \eqref{eq:ExteriorHarnackPlanarGeneral} for $x \in D(0,2R)^c$. Decompose the hitting distribution on whether or not we enter $D(0,2R)$ via the $R$-annulus: uniformly for $x \in D(0,2R)^c$ and $y \in \bd D(0,r)_{s}$, 
\begin{align*}
H_{D(0,r+s)}(x,y) & = P^x( S_{T_{D(0,r+s)}} = y \, , \, T_{\bd D(0,R)_{R}} > T_{D(0,R)}) \\
 & \,\, + P^x( S_{T_{D(0,r+s)}} = y \, , \, T_{\bd D(0,R)_{R}} < T_{D(0,R)}).
\end{align*}
We can bound the first term by \eqref{eq:FP(2.71)}: 
\begin{align*}
P^x( S_{T_{D(0,r+s)}} = y \, , \, T_{\bd D(0,R)_{R}} > T_{D(0,R)}) & \leq P^x(T_{\bd D(0,R)_{R}} > T_{D(0,R)}) \\
 \leq \, c R^2 (\log R)^2 (R^{-M} + R^{-M}) & \leq c R^{-M + 2} (\log R)^2 < cr^{-M+2} (\log r)^2.
\end{align*}
By the strong Markov property at $T_{\bd D(0,R)_{R}}$, the second term can be bounded, uniformly for $x' \in \bd D(0,R)_{R}$, by \eqref{eq:ExteriorHarnackPlanarGeneral}: 
\begin{align*}
P^x &( S_{T_{D(0,r+s)}} = y \, , \, T_{\bd D(0,R)_{R}}< T_{D(0,R)}) \\
 & = \EV^x( H_{D(0,r+s)}(S_{T_{\bd D(0,R)_{R}}}, y), \, T_{\bd D(0,R)_{R}}< T_{D(0,R)}) \\
 & \leq (1 + O(m^{-1} \log m)) H_{D(0,r+s)}(x', y).
\end{align*}
Thus, combining the two, we have for $x \in D(0,2R)^c$ and $x' \in \bd D(0,R)_{R}$, 
\begin{align} 
H_{D(0,r+s)}(x,y) & = \left(1 + O\left(m^{-1} \log m\right)\right) H_{D(0,r+s)}(x',y) + O(r^{-M+2} (\log r)^2), \label{eq:ExteriorHarnackPlanarFarOutPrebound}
\end{align}
which gives \eqref{eq:ExteriorHarnackPlanarGeneral} for $x \in D(0,2R)^c$ and $x' \in \bd D(0,R)_R$. 
Applying \eqref{eq:ExteriorHarnackPlanarGeneral} again for the same $x'$ and $x'' \in D(0,2R)^c$ gives \eqref{eq:ExteriorHarnackPlanarGeneral} for $x, x'' \in D(0,2R)^c$.

To prove \eqref{eq:ExteriorHarnackPlanarZoomOutGeneral} for $x, x' \in \bd D(0,R)_{\sqrt{R}}$, decompose $H_{D(0,r+s)}(x, y)$ over the event $\{T_{D(0,r+s)} > T_{D(0,4mR)^c}\}$ 
%$\{T_{D(0,r+s)} > T_{D(0,4mR)^c}, T_{D(0,r+s)} = T_{\bd D(0,r)_s}\}$ 
to get 
\begin{align}
P^x&( S_{T_{D(0,r+s)}} = y \, , \, T_{D(0,r+s)} < T_{D(0,4mR)^c}) \label{eq:ExtHarnDecomp1} \\
%P^x&( S_{T_{D(0,r+s)}} = y \, , \, T_{D(0,r+s)} < T_{D(0,4mR)^c}, \, T_{D(0,r+s)} = T_{\bd D(0,r)_s}) 
 & = H_{D(0,r+s)}(x,y) - P^x( S_{T_{D(0,r+s)}} = y \, , \, T_{D(0,r+s)} > T_{D(0,4mR)^c}). \notag
\end{align}
By the strong Markov property at $T_{D(0,4mR)^c}$ and \eqref{eq:ExteriorHarnackPlanarGeneral}, the last term of \eqref{eq:ExtHarnDecomp1} can be further decomposed to 
\begin{align}
P^x&( S_{T_{D(0,r+s)}} = y \, , \, T_{D(0,r+s)} > T_{D(0,4mR)^c}) \label{eq:ExtHarnDecomp2} \\
 & = \EV^x( H_{D(0,r+s)}(S_{T_{D(0,4mR)^c}} , y) \, ; \, T_{D(0,r+s)} > T_{D(0,4mR)^c}) \notag\\
 & = (1 + O(m^{-1} \log m)) H_{D(0,r+s)}(x,y) P^x(T_{D(0,r+s)} > T_{D(0,4mR)^c}), \notag  
\end{align}
which gives us the first equality in \eqref{eq:ExteriorHarnackPlanarZoomOutGeneral}. 
%The second follows from \eqref{eq:ExteriorHarnackPlanarGeneral} and \eqref{eq:FP(2.21)}, since, for $x, x' \in \bd D(0,R)_{s}$, the maximum error on  \eqref{eq:FP(2.21)} is, if $|x| = R$ and $|x'| = R+s$, 
%\begin{align*}
%P^{x'}&(S_{T_{D(0,r+s)}} = y \, , \, T_{D(0,r+s)} < T_{D(0,4mR)^c}) - P^x(S_{T_{D(0,r+s)}} = y \, , \, T_{D(0,r+s)} < T_{D(0,4mR)^c}) \notag\\
% & \leq P^{x'}(T_{D(0,r+s)} < T_{D(0,4mR)^c}) - P^x(T_{D(0,r+s)} < T_{D(0,4mR)^c}) \notag\\
% & = \frac{\log\left(\frac{R+s}{r+s}\right) - \log\left(\frac{R}{r+s}\right) + O(r^{-1/4})}{\log\left(\frac{4mR}{r+s}\right)}
% = \frac{\log\left(\frac{R+s}{R}\right) + O(r^{-1/4})}{\log\left(\frac{4mR}{r+s}\right)} \notag\\
% & \leq cr^{-1/4}(\log m)^{-1} = o(m^{-1} \log m). \qed \notag
%\end{align*}
The second follows from \eqref{eq:FP(2.21)} and \eqref{eq:logEstimate}, since, for $x, x' \in \bd D(0,R)_{\sqrt{R}}$, if $|x| = R$ and $|x'| = R+\sqrt{R}$, 
\begin{align}
&\frac{P^{x'}(T_{D(0,r+s)} < T_{D(0,4mR)^c})}{P^x(T_{D(0,r+s)} < T_{D(0,4mR)^c})} 
 = \frac{\log\left(\frac{R+\sqrt{R}}{r+s}\right) + O(r^{-1/4})}{\log\left(\frac{R}{r+s}\right) + O(r^{-1/4})} \label{eq:EntranceBounds} \\
 & \,\,\, = 1 + O\left(\frac{\sqrt{R}}{R\log(\frac{R}{r+s})}\right) + o(r^{-1/4}) = 1 + o(m^{-1} \log m). \qed \notag
\end{align}

\newcommand{\Drs}{D(0,r+s)}
\newcommand{\hDrs}{\hp(D(0,r+s))}
\newcommand{\TDrs}{T_{D(0,r+s)}}
\newcommand{\hTDrs}{T_{\hp(D(0,r+s))}}

When attempting to move the planar exterior Harnack inequality to the torus, we run into difficulties in dealing with walks that wander and enter far-off copies of $D(0,r+s)$ instead of the primary copy. We modify the exterior Harnack inequality for the toral case to fit our requirements.

\begin{prop} \label{lem:ExteriorHarnackToral}
Let $R = 4mr$ with $1 \leq m = o(r^{1/4})$ and large enough $r$, $4mR < K/4$, and $s \leq (\log R)^4$. Then, uniformly for $\hat x, \hat x' \in \hp(\bd D(0,R)_{\sqrt{R}})$ and $\hat y \in \hp(\bd D(0,r)_{s})$, 
\begin{align}
P^{\hat x}& (\hat{S}_{\hTDrs} = \hat y; \, \hTDrs < T_{{\hp(D(0,4mR)}^c_K)}) \label{eq:ExteriorHarnackToralZoomOutGeneral} \\
&  = \left(1 + O\left(m^{-1} \log m\right)\right) P^{\hat x'}(\hat{S}_{T_{\hp(D(0,r+s))}} = \hat y; T_{\hp(D(0,r+s))} < T_{{\hp(D(0,4mR)}^c_K)}). \notag
\end{align}
\end{prop}

\pf For brevity, set 
\begin{align*}
D^* & := \Drs \cup D(0,4mR)^c, & & A_p := \{S_{T_{D^*}} = y\},   \\
\hat{D}^* & := \hp(\Drs) \cup \hp(D(0,4mR)^c_K),  & & A_t := \{\hat{S}_{T_{\hat{D}^*}} = \hat y\}.
\end{align*}
We start our walk at the primary copy $x$, consider the planar landing at the primary copy $y$, and decompose $P^x(A_t) = \hat{H}_{\hat{D}^*}(x, y)$ along the planar large disc escape time $T_{D(0,4mR)^c}$ and the toral annulus escape time $T_{\hat{D}^*}$:
\begin{align}
P^x(A_t) & = P^x(A_t; \, T_{\hat{D}^*} < T_{D(0,4mR)^c}) + P^x(A_t; \, T_{\hat{D}^*} \geq T_{D(0,4mR)^c}). \label{eq:ToralExtDecomp}
\end{align}
Since $\hp\inv\hat{D}^* \subset D^*$, $T_{D^*} \leq T_{\hat{D}^*}$ a.s. The first term of \eqref{eq:ToralExtDecomp} happens in the event $\{ T_{D^*} = T_{\hat{D}^*} = T_{\Drs} \}$, so the entirety of its action before the final step is inside the primary copy of $D(0,4mR)$. Hence, 
\[  P^x(A_t; \, T_{\hat{D}^*} < T_{D(0,4mR)^c}) = P^x(A_t; \,  T_{D^*} = T_{\hat{D}^*}) = P^x(A_p). \]
Note that $P^x(A_p)$ is \eqref{eq:ExteriorHarnackPlanarZoomOutGeneral}. The second term of \eqref{eq:ToralExtDecomp} only occurs if a targeted jump lands in a non-primary copy of $D(0,4mR) \setminus D(0,r+s)$. Hence, by \eqref{eqn:ProbTorusExit}, 
\begin{align*}
P^x(A_t; \, T_{\hat{D}^*} \geq T_{D(0,4mR)^c}) & \leq P^{x}(T_{D(0,4mR)^c} < T_{\hp(D(0,4mR)^c_K)}) \\
 & \leq O(K^{-M} (mR)^2). 
\end{align*}
\eqref{eq:ToralExtDecomp} thus reduces to $P^x(A_t) = P^x(A_p) + O(K^{-M} (mR)^2)$, which, by \eqref{eq:ExteriorHarnackPlanarZoomOutGeneral}, is 
\begin{align}
P^x(A_t) = (1 + O(m^{-1} \log m)) P^{x'}(A_t) + O(K^{-M} (mR)^2). \label{eq:AtApExterior}
\end{align}
Since $M > 4$, the error term $O(K^{-M} (mR)^2) = o(K^{-2-\beta})$ is absorbed via \eqref{eq:ExteriorErrorBound} applied to the $P^x(A_p)$ term above, with \eqref{eq:HalfwayFractalGambler} one ``level'' up ($D(0,4mR)^c$ as the outer bound instead of $D(0,R)^c$, $D(0,r+s)$ instead of $D(0,\frac{r}{4m}+s)$, %\njc {\bf  if you change this above to $\frac{r}{4m}$ change it here also}, 
and $x, x' \in \bd D(0,R)_{\sqrt{R}}$ instead of $\bd D(0,r)_r$), which yields \eqref{eq:ExteriorHarnackToralZoomOutGeneral}. $\qed$

\singlespacing
%\include{Bibliography}
%\cleardoublepage

\end{document}